\documentclass[twoside]{article}
\usepackage{Latex-document}
\usepackage{amsmath,amscd,amssymb,amsbsy,euscript}
\def\text{\mbox}

\def\bracket{[]}
\def\brace{\cdot\{\}}
\def\phom{\text{\bf prohom}}
\def\pro{{\bf pro}\;}
\def\red{\text{\bf red}}

\def \d{{\frak d}}
\def\defo{\text{\bf def}}

\def \g {{\frak g}}
\def \m{{\frak m}}
\def\n{{\frak n}}
\def \G{{\cal G}}
\def\P{{\cal P}}
\def\h {{\frak h}}
\def\exctimes{\stackrel !\otimes}
\def\hcup{{\cup_h}}

\def \D{{\cal D}}

\def\liminv{\text{liminv}\;}
\def\limdir{\text{limdir}\;}

\def \A{{\cal A}}

\def \Hom{\hom}

\def\M{{\cal M}}
\def\Lie{\text{\bf lie}}
\def \lie{
\Lie}
\def \O{{\cal O}}
\def \ch{\text{\bf ch}}

\newtheorem{Theorem}{Theorem}[section]

\newtheorem{Conjecture}[Theorem]{Conjecture}

\newtheorem{Proposition}[Theorem]{Proposition}

\newtheorem{Lemma}[Theorem]{Lemma}

\renewcommand{\thesection}{\arabic{section}}
\renewcommand{\thesubsection}{\thesection.\arabic{subsection}.}
\def\addsec{\addtocounter{section}{1} \setcounter{Theorem}{0}}

\title{\bf Deformations of Chiral Algebras\vskip 6mm}
\author{Dimitri Tamarkin\vspace*{-0.5cm}\thanks{Department of Mathematics, Harvard University,
1 Oxford Street, Cambridge, MA 02138, USA. E-mail: tamarkin@math.harvard.edu}}
\date{\vspace{-8mm}}

\begin{document}

\maketitle

\thispagestyle{first} \setcounter{page}{105}

\begin{abstract}

\vskip 3mm

We start studying chiral algebras (as defined by A. Beilinson and V. Drinfeld) from the point of view of
deformation theory. First, we define the notion of deformation of a  chiral algebra on a smooth curve $X$ over a
bundle of local artinian commutative algebras on $X$ equipped with a flat connection (whereas `usual' algebraic
structures are deformed over a local artinian algebra) and we show that such deformations are controlled by a
certain *-Lie algebra $\mathfrak g$. Then we try to contemplate  a possible additional structure on $\mathfrak g$
and we conjecture that this structure up to homotopy is a chiral analogue of Gerstenhaber algebra, i.e. a coisson
algebra with odd coisson bracket (in the terminology of Beilinson-Drinfeld). Finally, we discuss possible
applications of this structure to the problem of quantization of coisson algebras.

\vskip 4.5mm


\end{abstract}


\section*{1. Introduction} \addsec

\vskip-5mm \hspace{5mm}

Chiral algebras were introduced in \cite{BD}.
In the same paper the authors introduced the classical limit of a chiral algebra
which they call a coisson algebra and posed the problem of quantization
of coisson algebras. The goal of this paper is to show how
the theory of deformation quantization (=the theory
of deformations of associative algebras of a certain type)
in the spirit of \cite{K}
 can be developed in this situation.

Central object in the theory of deformations of  associative algebras
is the differential graded Lie algebra of Hochschild cochains.
It turns out that in our situation
it is more appropriate to use what we call pro-*-Lie-algebras rather than usual Lie algebras
(the notion of *-Lie algebra
was also introduced in \cite{BD}).
Next, we compute the cohomology of the pro-*-Lie-algebra controlling chiral deformations
of a free commutative $\D_X$-algebra $SK$, where $K$ is a locally free $\D_X$-module.

Next,  we state an analogue of  Gerstenhaber theorem which says
that the cohomology of the deformation complex of an associative algebra
carries the structure of a Gerstenhaber algebra. We give a definition of
a chiral analogue of Gerstenhaber algebra and define the operations of this structure
on deformation pro-*-Lie algebra of a chiral algebra.

Finally, mimicking  Kontsevich's formality theorem, we formulate
the formality conjecture for the deformation pro-*-Lie algebra of
the chiral algebra $SK$ mentioned above and claim that
this conjecture implies a 1-1 correspondence between
deformations of $SK$ and coisson brackets on $SK$.

\section*{2. Chiral algebras and their deformations} \addsec

\vskip-5mm \hspace{5mm}

\subsection{Chiral operations}

\vskip-5mm \hspace{5mm}

In \cite{BD} chiral  operations are defined as follows.
Let $X$ be a smooth curve and $M_i$,$N$ $\D_X$-modules.
Denote by $i_n:X\to X^{n}$ the diagonal embedding and by
$j_n:U_n\to X^n$ the open embedding of the complement to all
diagonals in $X^n$.
Set
\begin{equation}\label{def}
P_{\ch}(M_1,\ldots, M_n;N)=\Hom_{\D_{X^n}}(j_*j^*(M_1\boxtimes\cdots\boxtimes M_n),i_{n*}N).
\end{equation}
In the case $n=0$ set $$ P_{\ch}(M)=H^0(M\otimes_{\D_X}\O_X).
$$
Let $M$ be a fixed $D_X$-module.
Write
$$
P_{\ch M}(n)=
P_{\ch}(M,M,\ldots,M;M).
$$
It is explained in \cite{BD} that  $P_M$ is an operad.

\subsubsection{Chiral algebras}

\vskip-5mm \hspace{5mm}

Let $\lie$ be the operad of Lie algebras.
{\em A chiral algebra} structure on $M$ is a homomorphism
$\lie\to P_M$.  We have a standard chiral algebra structure on
$M=\omega_X$. A chiral algebra $M$ is called unital if it is endowed
with an injection $\omega_X\to M$ of chiral algebras.

\subsection{Deformations}

\vskip-5mm \hspace{5mm}

\subsubsection{Agreements} \label{agreements}

\vskip-5mm \hspace{5mm}

To simplify the exposition, we will only consider unital chiral algebras $M$ with the following restrictions: we assume
that $X$ is affine and the $D_X$-module  $M$ can be represented
 as $M=\omega_X\oplus N$, where $N\cong E\otimes_{\O_X} D_X$
for some locally free coherent sheaf $E$.

\subsubsection{Nilpotent \boldmath{$\D_X$}-algebras}

\vskip-5mm \hspace{5mm}

Let $E$ be a  left
$\D_X$-module  equipped with  a  commutative associative unital
product $E\otimes E\to E$.
 Let
$u:\O_X\to E$ be the unit embedding.
Call $E$ nilpotent if there exists a  $\D_X$-module splitting
$s:E\cong \M\oplus \O_X$ and  a positive integer $N$ such that the $N$-fold
product vanishes on $\M$. $\M$ is then a unique maximal $\D_X$-ideal
in $E$.

\subsubsection{Deformations over a  nilpotent  \boldmath{$\D_X$}-algebra}

\vskip-5mm \hspace{5mm}

Let  $E$ be
a nilpotent $\D_X$-algebra
with maximal ideal $\M$.     We have a notion of $E$-module and of
an $E$-linear chiral algebra .
For any $\D_X$-module $M$, $M_E:=M\otimes_{\O_X} E$ is an $E$-module.

Let $M$ be  a chiral algebra.
An $E$-linear unital chiral   algebra structure on  $M_E$ is called  {\em deformation
of $M$ over $E$}
if the induced structure on $M_E/ \M.M_E\cong M$ coincides with the  one
on $M$.
Denote by $G_M(E)$ the set of all isomorphism classes of such deformations.

\subsection{The functor \boldmath{$G_M$} and its representability}

\vskip-5mm \hspace{5mm}

 It is clear that $E\mapsto G_M(E)$ is a functor from the category
of nilpotent $D_X$-algebras to the category of sets.
In classical deformation theory one usually has a functor from the category
of (usual) local Arminian (=nilpotent and finitely dimensional) algebras to the category of sets and one tries to
represent it by a differential graded  Lie algebra.  In this section we will see
that in our situation a natural substitute for a  Lie algebra  is a so-called
*-Lie algebra in the sense of \cite{BD}.
More precise, given a *-Lie algebra $\g$, we are going to construct a functor
$F_{\g}$ from the category of nilpotent $\D_X$-algebras to the category of
sets.
In the next section we will show that the
functor $G_M$ is 'pro-representable' in this sense. We will
construct a pro-*-Lie algebra $\defo_M$ (exact meaning will be given below)
 and an isomorphism of functors
$G_M$ and $F_{\defo_M}$.

\subsubsection{*-Lie algebras}

\vskip-5mm \hspace{5mm}

\cite{BD} Let $\g_i$, $N$ be  right $\D_X$-modules. Set
$$
P_*(\g_1,\ldots, \g_n;N):=\hom_{\D_{X^n}}(\g_1\boxtimes\cdots
\boxtimes \g_n,i_{n*}N),
$$
and
$P_{*\g}(n):=P(\g,\ldots,\g;\g)$.  It is known that $P_{*\g}$ is an operad.
A *-Lie algebra structure on $\g$ is by definition a morphism of operads
$f:\lie\to P_{*\g}$.  Let  $b\in
\lie(2)$ be the element corresponding to the  Lie bracket.
We call $f(b)\in P_{*\g}(2)$ the *-Lie bracket.

\subsubsection{}

\vskip-5mm \hspace{5mm}

Let $\g$ be a *-Lie algebra and $A$ be a commutative $\D_X$-algebra.
introduce a vector space $\g(A)=\g\otimes_{\D_X} A$.
This space is naturally a Lie algebra. Indeed, we have a *-Lie bracket
$\g\boxtimes \g\to i_{2*}\g.
$
Multiply both parts by $A\boxtimes A$:
 $$(\g\boxtimes \g)\otimes_{\D_{X\times X}}(A\boxtimes A)\to i_{2*}\g
\otimes_{\D_{X\times X}}(A\boxtimes A).\eqno(*)
$$
The left hand side is isomorphic to $\g(A)\otimes \g(A)$.
The right hand side is isomorphic
to $\g\otimes_{\D_X} (A\otimes_{\O_X} A)$.
Thus, (*) becomes:
$$\g(A)\otimes \g(A)\to \g\otimes_{\D_X} (A\otimes_{\O_X} A).
$$
The product on $A$
gives rise to a map
$$
\g\otimes_{\D_X} (A\otimes_{\O_X} A)\to \g\otimes_{\D_X} A\cong \g(A),
$$
and we have a map $\g(A)\otimes\g(A)\to\g(A)$. It is straightforward
to check that this map is a Lie bracket.

\subsubsection{}

\vskip-5mm \hspace{5mm}

Let now $\g$ be a differential  graded *-Lie algebra
and let  $A$ be a differential graded commutative $\D_X$-algebra. Then
$\g(A):=\g\otimes_{\D_X} A$ is a differential graded Lie algebra.

\subsubsection{}

\vskip-5mm \hspace{5mm}

Let $A$ be a nilpotent $\D_X$ algebra and $\M_A\subset A$
be the maximal nilpotent ideal. Then $\g(\M_A) $ is  a nilpotent differential
graded Lie algebra.

\subsubsection{}

\vskip-5mm \hspace{5mm}

Recall that given a differential graded
nilpotent Lie algebra $\n$,
one can construct the so called Deligne groupoid $\G_\n$.
 Its objects are all $x\in \n^1$ satisfying $dx+[x,x]/2=0$
(so called Maurer-Cartan elements).  The group
$\text{exp}(\n^0)$  acts on the set of Maurer-Cartan elements
by gauge transformations. $\G_\n$ is the groupoid of this action.
Denote by $\D_\n$ the set of isomorphism classes of this groupoid.
If $f:\n\to \m$ is a map of differential graded Lie algebras
such that the induced map on cohomology $H^i(f)$ is an isomorphism
for all $i\geq 0$,
then the induced map $\D_\n\to \D_\m$ is a bijection.
If $\n,\m$ are both centered in non-negative degrees, then
the induced map $\G_\n\to\G_\m$ is an equivalence of categories.
Since in our situation we will deal with Lie algebras centered in arbitrary
degrees, we will   use $\D_\n$ rather than  groupoids.

\subsubsection{}

\vskip-5mm \hspace{5mm}

Set $F_M(A)=\D_{\g(\M_A)}$. It is a functor from the category of
nilpotent $\D_X$-algebras to the category of sets.

\subsection{Pro-*-Lie- algebras}

\vskip-5mm \hspace{5mm}

 *-Lie algebras are insufficient for description of deformations of
chiral algebras. We will thus develop  a generalization.
We need some preparation

\subsubsection{Procategory}

\vskip-5mm \hspace{5mm}

For an Abelian category $C$ consider the category $\pro C$ whose
objects are functors $I\to C$, where $I$ is a small filtered
category. Let $F_k:I_k\to C$, $k=1,2$ be objects.
Set $$
\Hom(F_1,F_2):=\liminv_{i_2\in I_2}\limdir_{i_1\in I_1} (F_1(i_1),F_2(i_2)).
$$
The composition of morphisms
is naturally defined. One can show that $\pro C$ is an Abelian category.
Objects of $\pro C$ are called pro-objects.

\subsubsection{Direct image of pro-$\D$-modules}

\vskip-5mm \hspace{5mm}

Let $M:I\to\D_Y-mod$
be a pro-object, where $Y$ is a smooth algebraic variety and let $f:Y\to Z$ be a locally closed embedding.
Denote the composition $f_*\circ M:I\to \D_Z-mod$   simply by
$f_*M$. We will get a functor $f_*:\pro\D_Y-mod\to\pro\D_Z-mod$.

\subsubsection{Chiral and *-operations}

\vskip-5mm \hspace{5mm}

For $N,M_i\in \pro\D_X-mod$ we define
$P_{*}(M_1,\ldots, M_n,N)$, $P_{\ch}(M_1,\ldots, M_n,N)$
by exactly the same formulas as for usual $\D_X$-modules.

\subsubsection{pro-*-Lie algebras}

\vskip-5mm \hspace{5mm}

*-Lie algebra structure on a pro-$\D_X$-module
is defined in the same way as for usual $\D_X$-modules.

\subsubsection{}

\vskip-5mm \hspace{5mm}

For a pro-right $\D_X$-module $I\to M$ and a left $\D_X$-module
$L$ define a vector space $M\otimes_{\D_X}L=\liminv_{I}(M\otimes_{\D_X} L)$.
For a *-Lie
algebra $\g$ and a commutative $\D_X$-algebra $a$,
 $\g\otimes_{\D_X} a$ is a Lie algebra. Construction is the same
as  for usual *-Lie algebras.
Similarly, we can define the functor
$F_{\g}$ from the category of nilpotent $\D_X$-algebras to the
category of sets.

\subsection{Representability of $G_M$ by a pro-*-Lie algebra}

\vskip-5mm \hspace{5mm}

We are going to construct a differential graded *-pro-Lie algebra $\g$ such that $F_\g$ is equivalent to $G_{\M}$.
We need a couple of technical lemmas.

\subsubsection{}

\vskip-5mm \hspace{5mm}

Let $Y$ be a smooth affine algebraic varieties and $U,V$ be right $\D_Y$-modules. Let $U_{\alpha},\alpha\in A$ be
the family of all finitely generated submodules of $U$. Denote  $\phom(U,V)=\liminv_{\alpha}(U_{\alpha},V)$ the
corresponding pro-vector space.

\subsubsection{}

\vskip-5mm \hspace{5mm}

Let $i:X\to Y$ be a
closed embedding, let $B$ be a $\D_Y$-module and
$M$ be a $\D_X$-module.
Then
$$
\phom_{\D_Y}(B,i_*(M\otimes_{\O_X}\D_X))
$$
is a pro-$\D_X$-module. Denote it by $P(B,M)$.
Let now $Y=X^n$.
\begin{Lemma}
Assume that $B=j_{n*}j_n^*(E\otimes_{\O_{X^n}}\D_{X^n})$,
where $E$ is locally free and coherent.
   For any left $\D_X$-module $L$ we have
$$
\phom(B,i_{n*}(M\otimes_{\O_X} L))\cong P(B,M)\otimes_{\D_X} L.
$$
\end{Lemma}

{\bf Proof.} Let $F=j_{n*}j_n^*E$. We have $B=F\otimes_{\O_{X^n}}{\D_{X^n}}$. Note that $F=\limdir F_\alpha$,
where $F_\alpha$ runs through the set of all free coherent submodules of $F$.

We
have
\begin{eqnarray*}
P(B,M)=\liminv\hom_{\D_{X^n}}(F_\alpha\otimes_{\O_{X^n}} \D_{X^n}, i_{n*}(M\otimes_{
\O_X} L))\\
\cong
\liminv F_\alpha^*\otimes_{\O_{X^n}} i_{n*}(M\otimes_{\O_X} \D_X)\otimes_{\D_X} L\\
\cong\liminv \hom_{\O_{X^n}}(F_\alpha,i_{n*}(M\otimes_{\O_X} \D_X))\otimes_{\D_X} L\\
\cong\phom(B,i_{n*}(M\otimes_{\O_X} \D_X))\otimes_{\D_X} L.
\end{eqnarray*}

\subsubsection{} \label{dirim}

\vskip-5mm \hspace{5mm}

Let $B,M$ be as above. We have a natural morphism
$$
i:i_{n*}P(B,M)\cong P(B,M)\otimes_{\D_X}\D_X^{\otimes_{\O_X}^n} \to \phom(B, M\otimes_{\D_X}\otimes
(\D_X)^{\otimes_{\O_X}^n}).
$$
The above lemmas imply that $i$ is an isomorphism.

\subsubsection{}

\vskip-5mm \hspace{5mm}

Let $M$ be a right $\D_X$-module.
 Set
$$U_M(n)=\P_{\ch}(M,M,\ldots,M; M\otimes \D_X):=
\phom(j_{n*}j_n^{*}M^{\boxtimes n},i_{n*}(M\otimes \D_X)),
$$
it is a right pro-$\D_X$-module.
We will endow  the collection $U_M$ with the structure  of an operad in *-pseudotensor category.
This means that we will define the composition maps
$$\circ_i\in P_*(U_M(n), U_M(m);U_M(n+m-1)),
$$
$i=1,\ldots,n+m-1$, satisfying the operadic axioms.
We need a couple of technical facts.

\subsubsection{}

\vskip-5mm \hspace{5mm}

Let $i_n:X\to X^n$ be the diagonal embedding
and $p_n^i:X^n\to X$ be the projections.
Lemma \ref{dirim} implies that
\begin{Lemma}
$$
i_{n*}U_M(k)\cong \P_{\ch}(M,\ldots,M; M\otimes_{\D_X} \D_X^{\otimes n}).
$$
\end{Lemma}
\begin{Lemma} \label{di}
For any $\D_X$-modules $M,S$ we have
 an isomorphism
$$
i_{n*}(M)\otimes p_n^{i*}S\cong i_{n*}(M\otimes S).
$$
\end{Lemma}

\subsubsection{}

\vskip-5mm \hspace{5mm}

We are now ready
to define the desired structure. In
virtue of \ref{di}
we have natural maps:
$$
\P_{\ch}(M_1,\ldots,M_n;N)\to
\P_{\ch}(M_1,\ldots, M_i\otimes \D_X,\ldots,M_n;(N\otimes \D_X)).
$$
Thus, we have maps:
\begin{eqnarray*}
& & \P_{\ch}(M_1,\ldots, M_n;(N_i\otimes \D_X))\boxtimes \P_{\ch}(N_1,\ldots N_m;(K\otimes \D_X))\\
& \to & \P_{\ch}(M_1,\ldots, M_n;(N_i\otimes \D_X)) \\
& & \boxtimes \P_{\ch}(N_1,\ldots, N_i\otimes \D_X,\ldots,
N_m;(K\otimes \D_X)\otimes \D_X) \\
& \to & \P_{\ch}(N_1,\ldots N_{i-1},M_1,\ldots,M_n,N_{i+1},\ldots, N_m;K\otimes \D_X\otimes \D_X)\\
& \cong & i_{2*}\P_{\ch}( N_1,\ldots N_{i-1},M_1,\ldots,M_n,N_{i+1},\ldots, N_m,K\otimes\D_X).
\end{eqnarray*}
By substituting $M$ instead of all $N_i,M_j,K$, we get the desired insertion map
$$
\circ_i:U_M(n)\boxtimes U_M(m)\to i_{2*}U_M({n+m-1}).
$$

\subsubsection{}

\vskip-5mm \hspace{5mm}

Similarly, we have insertion maps
$$
\circ_i:U_M(n)\otimes \P_{\ch M}(m)\to U_M(n+m-1),
$$
and
$$
\circ_i:P_{\ch M}(n)\otimes U_M(m)\to U_M(n+m-1).
$$

\subsubsection{}

\vskip-5mm \hspace{5mm}

Let $\O$ be a differential graded operad. Set
$$\g_{\O,n}:=\O(n)^{S_n},
$$
and $\g_{\O}=\oplus_n \g_{\O,n}[1-n]$.

Let $p_n:\O(n)\to\g_{\O,n}$ be the natural projection, which is the symmetrization map.
Define the brace $(x,y)\mapsto x\{y\}$,
$\g_{\O,n}\otimes \g_{\O,m}\to \g_{\O,n+m-1}$ as follows.
\begin{equation}\label{brace}
x\{y\}=np_{n}(\circ_1(x,y))
\end{equation}
and the bracket
\begin{equation}\label{lie}
[x,y]=x\{y\}-(-1)^{|x||y|}y\{x\}.
\end{equation}
We see that $[,]$ is a Lie bracket. Thus, $\g_\O$ is a differential graded Lie
algebra.
For an operad $\O$ denote by $\O'$ the shifted operad such that
the structure of an $\O'$-algebra on a complex  $V$ is equivalent
to the structure of an $\O$-algebra on a complex $V[1]$.
Thus, $\O'(n)=\O(n)\otimes \epsilon_n[1-n]$, where
$\epsilon_n$ is the sign representation of $S_n$.

Let $\O$ be an operad of vector spaces.
The set of Maurer-Cartan elements of $\g_{\O'}$ is in 1-1 correspondence with
maps of operads $\lie\to \O$.

Assume that $O(1)$ is a nilpotent algebra ($x^n=0$ for any $x\in \O(1)$).
Let $A$ be $O(1)$ with adjoined unit and let $A^\times $ be the group of invertible
elements. $A^\times$ acts on $\O$ by automorphisms. Therefore,
$A^\times$ acts on the set of maps $\lie\to\O$. The groupoid of this action
is isomorphic to the Deligne groupoid of $\g_{\O'}$.

\subsubsection{}

\vskip-5mm \hspace{5mm}

Similarly, let $\A$ be a *-operad. Then  formula
\ref{lie} defines a Lie-* algebra $\g_{\A}$.  We have natural action of a usual
pro-Lie algebra $\g_{\P_{\ch}(M)}$ on a pro *-Lie algebra  $ \g_{U_M}$  by derivations.
The chiral bracket $b\in \g_{\P_{\ch}(M)}^1$ satisfies $[b,b]=0$.
Therefore, the bracket with $b$ defines a differential on
$\g_{U_M}$. Denote this differential graded *-Lie algebra by $\d_M$.

\subsubsection{}

\vskip-5mm \hspace{5mm}

To avoid using derived functors, we will slightly modify $\d_M$.
Recall that $M=\omega_X\oplus N$, where $N$ is free. Let
$$
\P_{\ch }^{\red}(M,\ldots,M;M\otimes \D_X)
\subset \P_{\ch}(M,\ldots,M;M\otimes \D_X)
$$ be the subset of all operations vanishing under all restrictions
$$
\P_{\ch }(M,\ldots,M;M\otimes \D_X)
\to \P_{\ch}(M,\ldots,M,\omega_X,M,\ldots,M;M\otimes \D_X).
$$

Let $\defo_M\subset \d_M$ be the submodule such that
$$
\defo_M=\oplus_n (\P_{\ch }^{\red}(M,\ldots,M;M\otimes \D_X)
\otimes \epsilon_n)^{S_n}[1-n].
$$
We see that $\defo_M$ is a *-Lie differential subalgebra of $\d_M$.

\subsubsection{}

\vskip-5mm \hspace{5mm}

\begin{Proposition} The functors $G_M$ and $F_{
\defo_M}$ are canonically isomorphic.
\end{Proposition}

\subsection{Example} \label{example}

\vskip-5mm \hspace{5mm}

Let $K$ be a free left $\D_X$-module.
Let $T^iK=K^{\otimes_{\O_{X^i}}}$. The symmetric group
$S_i$ acts on the $\D_X$-module $T^iK$; let
$S^iK=(T^iK)^{S_i}$ be the submodule of invariants
and $SK=\oplus_{i=0}^{\infty}S^iK$. $SK$ is naturally a free
commutative $\D_X$-algebra and, hence,  $SK^r:=SK\otimes \omega_X$
is a chiral algebra. We will compute the cohomology of the $\D_X$-module
$\defo_{SK^r}$. Let $S_0K=\oplus_{n=1}^\infty S^iK$. We have:
$$
\defo_{SK^r}=\oplus_n (P_{\ch}(S_0K^r[1],\ldots,S_0K^r[1];SK^r\otimes \D_X)[1])^{S_n}.
$$

On the other hand, denote by $\Omega:=SK\otimes K$. Consider
$\Omega$ as an $SK$-$\D_X$-module of differentials of
$SK$. We have the de Rham
differential $D:S_0K\to \Omega$.
We have a through map
\begin{eqnarray*}
c_n:P_{\ch}(K[1]^r,\ldots,K[1]^r, SK[1]^r)^{S_n}\cong
P_{\ch}^{SK}(\Omega[1]^r,\ldots,\Omega[1]^r,SK[1]^r)^{S_n}\\
\stackrel D\to
P_{\ch}(S_0K[1]^r,\ldots,S_0K[1]^r,SK[1]^r)^{S_n},
\end{eqnarray*}
where $P_{\ch}^{SK}$ stands for $SK$-linear chiral operations.
Denote by the same letter the induced
map
$$
c_n: P_{\ch}(K[1]^r,\ldots,K[1]^r, SK[1]^r)^{S_n}\to \defo_{SK}.
$$
\begin{Proposition}

\begin{enumerate}
\item[(1)]
$dc_n=0$;
\item[(2)]
$c_n$ induces an isomorphism
$$
P_{\ch}(K[1]^r,\ldots,K[1]^r, SK[1]^r)^{S_n}\to H^{n-1}(\defo_{SK})[1-n].
$$
\end{enumerate}
\end{Proposition}

\subsubsection{}

\vskip-5mm \hspace{5mm}

For a chiral algebra $M$ denote by $H_M$
the graded Lie algebra of cohomology of $\defo_M$.

\subsubsection{}

\vskip-5mm \hspace{5mm}

Assume that $K$ is finitely generated. Let
$$
K^\lor=\hom(K,\D_X)\otimes (\omega_X)^{-1}
$$
be the dual  module.
Then
$$
H_{SK^r}\cong \oplus_n
 (P_* (K^r,\ldots,K^r;SK^r\otimes \D_X)\varepsilon_n)^{S_n}[1-n]\\
=\oplus_n(\wedge^n K^\lor\otimes_{\O_X} SK)^r[1-n].
$$

\subsubsection{}

\vskip-5mm \hspace{5mm}

We will postpone the calculation of the *-Lie bracket on $H_{SK^r}$ until we show
in the next section that $H_M$ has in fact a richer structure.

\section*{3. Algebraic structure on the cohomology of the deformation
pro-*-Lie algebra} \addsec\setcounter{subsection}{0}

\vskip-5mm \hspace{5mm}

We will keep the agreements  and the notations from \ref{agreements}.

\subsection{Cup product}

\vskip-5mm \hspace{5mm}

We will  define
a chiral operation $\cup\in P_{\ch\; \defo_M[-1]}(2)$
and then we will study the induced map on cohomology.

\subsubsection{}

\vskip-5mm \hspace{5mm}

Recall that
$$
\defo_M[-1]\cong \oplus_n  (a_n)^{S_n},
$$
where
$$
a_n=P_{\ch}(N[1],\ldots,N[1];M\otimes \D_X).
$$
Let $i_n:X\to X^n$ be the diagonal embedding and let
$p^i:X^n\to  X$
$U_n\subset X^n$ be the complement to the union of all pairwise diagonals
$p^ix=p^jx$
and $j_n:U_n\to X_n$ be the open embedding.
Let $U_{n,m}\subset X^{n+m}$ be the complement to the diagonals
$p^ix=p^jx$, where $1\leq i\leq n$, $n+1\leq j\leq n+m$
and $j_{nm}:U_{nm}\to X^{n+m}$ be the embedding

Compute
\begin{eqnarray*}
j_{2*}j_2^*(a_n\boxtimes a_m)\cong
\hom(j_{n*}j_n^*(N[1]^{\boxtimes n})\boxtimes j_{m*}j_m^*(N[1]^{\boxtimes m}),\\
(i_{n}\times i_m)_*(j_{2*}j_2^*(M\otimes\D_X\boxtimes M\otimes \D_X)))\\
\cong \hom(j_{n*}j_n^*(N[1]^{\boxtimes n})\boxtimes j_{m*}j_m^*(N[1]^{\boxtimes m})
\otimes j_*\O(U_{n,m}),\\
(i_{n}\times i_m)_*(j_{2*}j_2^*(M\otimes\D_X\boxtimes M\otimes \D_X)\otimes j_*\O(U_{n,m})))\\
\cong\hom(j_{n+m*}j_{n+m}^*(N[1])^{\boxtimes n+m},(i_{n}\times i_m)_*(j_{2*}j_2^*(M\boxtimes M))\otimes \D_{X\times X})).\\
\end{eqnarray*}

Taking the composition with the chiral operation on $M$, we obtain a chiral operation
$$
j_{2*}j_2^*(a_n\boxtimes a_m)\to \hom(j_{n+m*}j_{n+m}^*(N[1])^{\boxtimes n+m},
i_{n+m*}(M\otimes \D_X\otimes \D_X))\cong i_{2*}a_{n+m},
$$
which induces a chiral operation from
$P_{\ch}(a_n^{S_n},a_m^{S_m};a_{n+m}^{S_{n+m}})$ and, hence, an  operation
$\cup\in (P_{\ch}(\defo_M[-1],\defo_M[-1];\defo_M[-1]))^{S_2}$.

\subsubsection{}

\vskip-5mm \hspace{5mm}

To investigate the properties of this operation, consider the brace *-operation
 $\brace\in P_{*}(\defo_M,\defo_M;\defo_M)$ defined by
formula (\ref{brace}).
Let $$r:P_{\ch}(A_1,A_2;A_3)\to P_{*}(A_1,A_2;A_3)$$
be the natural map
\begin{Proposition}
$d(\brace)=r(\cup)$.
\end{Proposition}

Let $\hcup$ be the induced operation on $H_M[-1]$. The above proposition implies that $r(\hcup)=0$. In virtue of
exact sequence
$$
0\to \hom((A_1)^l\otimes (A_2)^l,(A_3)^l)\to  P_{\ch}(A_1,A_2;A_3)\stackrel r\to P_*(A_1,A_2;A_3),
$$
$\hcup$ defines a $\D_X$-commutative product $H_M[-1]\otimes H_M[-1]\to H_M[-1]$, denoted by the same
letter.

\subsubsection{}

\vskip-5mm

\begin{Proposition}
$\hcup$ is associative.
\end{Proposition}

\subsubsection{Leibnitz rule}

\vskip-5mm \hspace{5mm}

We are going to establish a relation between $\cup$ and $\brace$.
This relation is similar to the one of coisson algebras.
Our exposition will mimic the definition of coisson algebras from  \cite{BD}.

\subsubsection{}

\vskip-5mm \hspace{5mm}

Let $A_i$ be right $\D_X$-modules.
Write $A_1\exctimes A_2:=(A_1^l\otimes A_2^l)^r$;
$P_!(A_1,A_2;A_3):=\hom(A_1\exctimes A_2,A_3)$.
We have
$(A\exctimes B)=i_2^*(B\boxtimes C)$;
$$
i_{2*}(A\exctimes B)\to i_{2*}A\otimes p_2^*(B^l).
$$

We have a map
$$
c:P_*(A_1,A_2;B)\otimes P_!(B,C;D)\to P_*(A_1,A_2\exctimes C;D)
$$ defined as follows. Let
$u:A_1\boxtimes A_2\to i_{2*}B$ and $m:B\exctimes C\to D$.
Put
$$
c(u,v):A_1\boxtimes (A_2\exctimes C)\cong (A_1\boxtimes A_2)\otimes p_2^*(C^l)\to
i_{2*}B\otimes p_2^*(C^l)\cong i_{2*}(B\exctimes C)\to i_{2*}D.
$$

\subsubsection{}\label{comp}

\vskip-5mm \hspace{5mm}

Denote
$$
e=c(\bracket,\hcup)\in P_*(H_M,H_M\exctimes H_M;H_M).
$$
Let $T:H_M\exctimes H_M\to H_M\exctimes H_M$
be the action of symmetric group and let
$e^T$ be the composition with $T$.
Let $f\in P_*(H_M,H_M\exctimes H_M;H_M)
$ be defined by:
$$
H_M\boxtimes(H_M\exctimes H_M)\stackrel\cup\to
H_M\boxtimes H_M\stackrel\bracket\to i_*H_M.
$$
\begin{Proposition}
 We have $f=e+e^T$.
\end{Proposition}

In other words, the cup product and the bracket satisfy the Leibnitz identity.

\subsubsection{}

\vskip-5mm \hspace{5mm}

We see that $H_M$ has a pro-*-Lie bracket,
$(H_M)^l[1]$ has a commutative $\D_X$-algebra structure,
and these structures satisfy the Leibnitz identity.
Call this structure {\it a c-Gerstenhaber algebra} structure.
 Thus, our findings can be summarized as follows.
\begin{Theorem}
The cohomology  of the deformation pro-*-Lie algebra of a chiral
algebra is naturally a pro-c-Gerstenhaber algebra.
\end{Theorem}

\subsection{Example $M=(SK)^r$}

\vskip-5mm \hspace{5mm}

We come back to our example \ref{example}. For simplicity assume $K$ is finitely
generated free $\D_X$-module. We have seen in this case that
$$
(H_M)^l[-1]\cong \oplus_i \wedge^i K^\lor\otimes S^K[-i]\cong S(K^\lor[-1]\oplus K).
$$
\begin{Proposition}
The cup product on $H_M$ coincides with  the natural one on the symmetric power
algebra.
\end{Proposition}

\subsubsection{}

\vskip-5mm \hspace{5mm}

To describe the bracket it suffices to define it on
the submodule of generators $G=(K^\lor[-1]\oplus K)^r$.
Define $\bracket\in P_*(G,G;H_M)$ to be zero when restricted onto
$K^r\boxtimes K^r$ and $K^{\lor r}[-1]\boxtimes K^{\lor r}[-1]$.
Restriction onto $K\boxtimes K^{\lor}[-1]$ takes values
in $\omega_X\subset H_M$ and is given by the canonic
$*$-pairing from \cite{BD}
$$
(K^{\lor}\boxtimes K)^r\to i_{2*}\omega_X.
$$
Recall the definition. We have $K^{\lor r}=\Hom(K^r,\D_X\otimes \omega_X)$. For open $U,V\subset X$
we have the composition map
$$
K(U)\otimes K^\lor(V)\to \D_X\otimes \omega_X(U\cap V)\cong
i_{2*}\omega_X(U\times V)
$$
which defines the pairing.
This uniquely defines the *-Lie bracket.

\section*{4. Formality Conjecture} \addsec\setcounter{subsection}{0}

\vskip-5mm \hspace{5mm}

Following the logic of Kontsevich's formality theorem, one can
formulate a formality conjecture in this situation.

\subsection{Quasi-isomorphisms}

\vskip-5mm \hspace{5mm}

A map  $f:\g\to\h$ of differential  graded pro-*-Lie algebras
is called quasi-isomorphism if it induces an isomorphism on cohomology.
Call a pro-*-Lie algebra {\em perfect} if such is its underlying complex of pro-vector spaces.
The morphism $f$ is called {\em perfect quasi-isomorphism} if it is a quasi-isomorphism and both
$\g$ and $\h$ are perfect.

Two  perfect pro-*-lie algebras are called perfectly quasi-isomorphic
if there exists a chain of perfect quasi-isomorphisms connecting $\g$ and $\h$.
\begin{Conjecture}
$\defo_{SK}$ and $H_{SK}$ are perfectly quasi-isomorphic.
\end{Conjecture}

The importance of this conjecture can be seen from the following theorem:

\begin{Theorem} Any chain of perfect quasi-isomorphisms between
$\defo_{SK^r}$ and $H_{SK^r}$ establishes a bijection between the set of
isomorphism classes of  $A$-linear coisson brackets
on $SK^r\otimes A$ which vanish modulo the maximal ideal $\M_A$
and the set of isomorphism classes
of all deformations of the chiral algebra $SK^r$ over $A$.
\end{Theorem}

\end{document}